\documentclass[12pt,bezier]{article}
\usepackage{times}
\usepackage{booktabs}
\usepackage{pifont}
\usepackage{floatrow}
\usepackage{caption}
\usepackage{mathrsfs}
\usepackage{amsmath}
\usepackage{amsfonts,amsthm,amssymb,mathrsfs,bbding}
\usepackage{txfonts}
\usepackage{graphics,multicol}
\usepackage{graphicx}
\usepackage{color}
\usepackage{amssymb}
\allowdisplaybreaks[3]
\usepackage{caption}
\usepackage{cite}
\usepackage{latexsym,bm}
\usepackage{indentfirst}
\usepackage{color}
\usepackage[colorlinks=true,anchorcolor=blue,filecolor=blue,linkcolor=blue,urlcolor=blue,citecolor=blue]{hyperref}
\usepackage{extarrows}
\usepackage{cite}
\usepackage{latexsym,bm}
\usepackage{mathtools}
\usepackage{enumerate}

\evensidemargin=\oddsidemargin\topmargin=-1.5cm

\newtheorem{thm}{Theorem}[section]

\newtheorem{lem}{Lemma}[section]

\theoremstyle{definition}

\addtocounter{section}{0}

\begin{document}
	\title{Sufficient conditions for the variation of toughness under the distance spectral in graphs involving minimum degree}
	\author{{\bf Peishan Li} \thanks{Corresponding author.
			E-mail addresses: lipsh2021@lzu.edu.cn (P.S. Li)}
		\\
		{\footnotesize  School of Mathematics and Statistics, Lanzhou University, Lanzhou, Gansu, China} \\}
	\date{}
	
	\date{}
	\maketitle
	\begin{abstract}
		 In 1988, Enomoto introduced a variation of toughness \(\tau(G)\) of 
			a connected non-complete graph $G=(V,E)$, defined as  
			\(\tau(G)=\min\left\{\frac{|S|}{c(G - S)-1}\right\}\), 
			in which the minimum  is taken over all proper sets $S\subset V(G)$, 
			where $c(G-S)$  denotes the number of components of $G$.  
		If \(\tau(G)\geq \tau\), then \(G\) is called a $\tau$-tough graph. 
		Chen, Fan and Lin (2025) provided sufficient spectral radius conditions 
		for a graph to be $\tau$-tough with minimum degree $\delta$. 
		Inspired by this, we propose  sufficient conditions for a graph 
		to be $\tau$-tough in terms of minimum degree 
		and distance spectral radius, where \(\tau\) or \(\frac{1}{\tau}\) is a positive integer.
		
		\medskip
		
		
		\medskip
		
	\end{abstract}

	\section{Introduction}
	
	Let $G=(V,E)$ be a graph with vertex set $V(G)$ and edge set $E(G)$.
		Its order is $|V(G)|$, denoted by $n$, and its size is $|E(G)|$,
		denoted by $m$. 
		We denote by $\delta(G)$ and $\Delta(G)$ the minimum degree 
		and the maximum degree of $G$, respectively. 
		Let $c(G)$ be the number of components of $G$. 
		For a subset $S\subseteq V(G)$, we
		use $G[S]$ and $G-S$ to denote the subgraphs of $G$ induced by $S$and $V(G)\backslash S$, respectively.  
		For two vertex-disjoint graphs $G_1$ and $G_2$, $G_1 \cup G_2$ 
		denotes the disjoint union of $G_1$ and $G_2$. 
		The join $G_1 \vee G_2$ is the graph obtained from $G_1 \cup G_2$ 
		by adding all possible edges between $V(G_1)$ and $V(G_2)$.
		
		The distance between vertices $v_i$ and $v_j$, denoted by $d_{ij}(G)$, 
		is the length of a shortest path from $v_i$ to $v_j$. 
		The distance matrix of $G$, denoted by $D(G)$, is the $n$-by-$n$ real symmetric matrix  
		whose $(i,j)$-entry is $d_{ij}(G)$. The distance eigenvalues in $G$ are the eigenvalues 
		of its distance matrix $D(G)$. 
		Let $\lambda_1(D(G)) \geq \lambda_2(D(G)) \geq \cdots \geq \lambda_n(D(G))$ 
		be its distance eigenvalues in nonincreasing order.  
		According to the Perron-Frobenius theorem, 
		$\lambda_1(D(G))$ is always positive 
		(unless $G$ is a trivial graph), and satisfies $\lambda_1(D(G)) \geq |\lambda_i(D(G))|$ for $i = 2, 3, \cdots, n$. We call $\lambda_1(D(G))$ the $distance$ $spectral$ $radius$ of $G$. 
	
	The concept of graph toughness was first introduced 
	by Chvátal \cite{Chvatal1973} in 1973. 
	For a connected non-complete graph $G$, the toughness $t(G)$ is defined as 
	$$t(G) = \min\limits_{S \subset V(G)} \left\{ \frac{|S|}{c(G-S)} : c(G-S) > 1 \right\}.$$ 
	If $t(G) \geq t$, then $G$ is called a $t$-tough graph. 
	In 1995, Alon first studied the relationship between the toughness 
	of $G$ and its eigenvalues. At the same time, Brouwer \cite{Brouwer1995} independently 
	proved a better lower bound $t(G)>\frac{d}{\lambda}-2$, 
	and conjectured that $t(G) \geq \frac{d}{\lambda}-1$, 
	where $\lambda$ is the second largest absolute eigenvalue 
	of the adjacency matrix. Later, Gu \cite{Gu2021a} proved part of this conjecture 
	and fully proved this conjecture\cite{Gu2021} in 2021.  
	Besides refining the lower bound of \(t(G)\), 
	it is interesting to extend the results on finding a condition 
	for a graph to be $t$-tough in terms of its  related eigenvalues. 
	Around this proposed, the researchers have begun to study the relationship 
	between the toughness of a graph and the spectral radius
	of other types of matrices, such as~\cite{Fan2023,Lou2025}.
	
	To better study the existence of factors in graphs, 
	Enomoto \cite{Enomoto1998} introduced a variation of toughness 
	\(\tau(G)\) of a graph $G$. 
	For a connected non-complete graph $G$, 
	$$\tau(G) = \min\limits_{S \subset V(G)}\left\{ \frac{|S|}{c(G-S)-1}: c(G-S) > 1 \right\}.$$ 
	If $\tau(G) \geq \tau$, then $G$ is called a $\tau$-tough graph. 
	It is natural to ask whether there are some sufficient conditions for a graph $G$ to be $\tau$-tough. 
	For this purpose, Chen, Fan and Lin \cite{Chen2024} (resp. Chen, Li and Xu~\cite{Chen2025}) 
	provided sufficient conditions based on the spectral radius (resp. signless Laplacian)
	for a graph $G$ to be $\tau$-tough, where 
	\(\tau\) or \(\frac{1}{\tau}\) is an positive integer. 
	Inspired by their work, we continue to explore the relationship 
	between $\tau$-tough and the distance spectral radius. 
	
	\begin{thm}\label{main1}
		Let $G$ be a connected graph of order 
			$n\geq \max\{9\delta, \frac{\delta^2}{2}+3\delta+3\}$ 
			with minimum degree $\delta\geq 2$. 
		If $\lambda_1(D(G)) \leq \lambda_1(D(K_\delta \vee (K_{n-2\delta-1} \cup (\delta+1)K_1)))$, then $G$ is a $1$-tough graph, 
		unless $G \cong K_\delta \vee (K_{n-2\delta-1} \cup (\delta+1)K_1)$.
	\end{thm}
	
	\begin{thm}\label{main2}
		Let $G$ be a connected graph of order $n\geq 4\tau^2+5\tau+1$, 
			where $\tau \geq 2$ be an integer. 
		If $\lambda_1(D(G)) \leq \lambda_1(D(K_{\tau-1} \vee (K_{n-\tau} \cup K_1)))$, then $G$ is a $\tau$-tough graph, unless $G \cong K_{\tau-1} \vee (K_{n-\tau} \cup K_1)$.
	\end{thm}
	
	\begin{thm}\label{main3}
	Let $G$ be a connected graph of order $n\geq 4\tau+\frac{1}{\tau}+5$, 
			where $\frac{1}{\tau} \geq 1$ be a positive integer. 
		If $\lambda_1(D(G)) \leq \lambda_1(D(K_1 \vee (K_{n-\frac{1}{\tau}-2} \cup (\frac{1}{\tau}+1)K_1)))$, then $G$ is a $\tau$-tough graph, unless $G \cong K_1 \vee (K_{n-\frac{1}{\tau}-2} \cup (\frac{1}{\tau}+1)K_1)$.
	\end{thm}

	\section{Preliminary lemmas}
In this section, we introduce some necessary notation, terminologies and useful results below.
		We begin with the relationship between the distance spectral radius of a graph and its spanning graph. 
	
	\begin{lem}[\cite{Godsil1993}]\label{le1}
		Let $e$ be an edge of a graph $G$ such that $G-e$ is connected. Then $$\lambda_1(D(G))<\lambda_1(D(G-e)).$$
	\end{lem}
	
	Let $M$ be a real $n\times n$ matrix. Suppose $M$ can be written in the following matrix
	\[
	M=\left(\begin{array}{ccccccc}
		M_{1,1}&M_{1,2}&\cdots &M_{1,m}\\
		M_{2,1}&M_{2,2}&\cdots &M_{2,m}\\
		\vdots& \vdots& \ddots& \vdots\\
		M_{m,1}&M_{m,2}&\cdots &M_{m,m}\\
	\end{array}\right),
	\]
	where the rows and columns are partitioned into subsets \( X_1, X_2, \cdots, X_m \) of \(\{1, 2, \cdots, n\}\). The quotient matrix \( R(M) \) of matrix \( M \) with respect to the given partition is an \( m \times m \) matrix. The \((i, j)\)-entry of \( R(M) \) is the average of the row-sums of the block \( M_{i,j} \) in \( M \). If the row-sum of each block \( M_{i,j} \) is a constant, then the corresponding partition is called an \emph{equivalent partition}, and the corresponding matrix is called an \emph{equivalent quotient matrix}.
	
	\begin{lem}[\cite{Brouwer2011,Godsil2001,Haemers1995}]\label{le2}
		Let $M$ be a real symmetric matrix and let $R(M)$ be its equitable quotient matrix.
		Then the eigenvalues of the quotient matrix $R(M)$ are eigenvalues of $M$.
		Furthermore, if $M$ is nonnegative and irreducible, then the spectral radius of the quotient matrix $R(M)$ equals to the spectral radius of $M$.
	\end{lem}
	
	For a connected graph \( G \) of order \( n \), let \( W(G) = \sum_{i<j} d_{ij}(G) \) be the Wiener index of graph \( G \). From the Rayleigh quotient \cite{Horn1985}, we obtain the following lemma.
	
	\begin{lem}\label{le5}
		Let $G$ be a connected graph with order $n$. Then
		\begin{equation*}
			\lambda_1(D(G))=\mathop{\max}_{\mathbf{x}\neq\mathbf{0}}\frac{\mathbf{x}^TD(G)\mathbf{x}}{\mathbf{x}^T\mathbf{x}}\geq \frac{\mathbf{1}^TD(G)\mathbf{1}}{\mathbf{1}^T\mathbf{1}}=\frac{2W(G)}{n},
		\end{equation*}
		where $\mathbf{1}=(1, 1, \ldots, 1)^T.$
	\end{lem}
	
	Let \(W^{(2)}(G)\) denote the sum of squared distances between all unordered vertex pairs in graph \( G \), that is,\[
	W^{(2)}(G)=\sum_{1\leq i<j\leq n}d_{ij}^2(G).
	\]
	
	\begin{lem}[\cite{Zhou2010,Zhou2007a,Zhou2007b}]\label{le6}
		Let $G$ be a graph on $n\geq 2$ vertices with sum of the squares of the distances between all unordered pairs of vertices $W^{(2)}(G).$ Then
		$$\lambda_1(D(G))\leq \sqrt{\frac{2(n-1)W^{(2)}(G)}{n}}$$
		with equality if and only if $G$ is the complete graph $K_n$, and if $G$ has exactly one positive distance eigenvalue, then
		$$\lambda_1(D(G))\geq \sqrt{W^{(2)}(G)}$$
		with equality if and only if $G$ is $K_2.$
	\end{lem}
	
	\begin{lem}[\cite{Zhang2021}]\label{le3}
		Let $n, c, s$ and $n_i~(1\leq i\leq c)$ be positive integers with $n_1\geq n_2\geq \cdots \geq n_c\geq 1$ and $n_1+n_2+\cdots +n_c=n-s.$ Then $$\lambda_1(D(K_s\vee(K_{n_1}+K_{n_2}+\cdots+K_{n_c})))\geq \lambda_1(D(K_s\vee(K_{n-s-(c-1)}+(c-1)K_1)))$$
		with equality if and only if $(n_1, n_2, \ldots, n_c)=(n-s-(c-1), 1, \ldots, 1).$
	\end{lem}
	
	\begin{lem}[\cite{Lou2025}]\label{le4}
		Let $n, c, s, p$ and $n_i~(1\leq i\leq c)$ be positive integers with $n_1\geq 2p$, $n_1\geq n_2\geq \cdots \geq n_c\geq p$ and $n_1+n_2+\cdots +n_c=n-s.$
		Then $$\lambda_1(D(K_s\vee(K_{n_1}+K_{n_2}+\cdots+K_{n_c})))\geq \lambda_1(D(K_s\vee(K_{n-s-p(c-1)}+(c-1)K_p)))$$
		with equality if and only if $(n_1, n_2, \ldots, n_c)=(n-s-p(c-1), p, \ldots, p).$
	\end{lem}
	
	\section{Proof of Theorem \ref{main1}}
	
	In this section, we give the proof of Theorem \ref{main1}.

	\medskip
	
	\noindent  \textbf{Proof of Theorem \ref{main1}.}
	Let $G^*=K_\delta \vee (K_{n-2\delta-1} \cup (\delta+1)K_1)$. 
	Suppose to the contrary that $G$ is not a $1$-tough graph. 
	By the definition of $1$-tough, 
	there exists a non-empty vertex subset $S \subset V(G)$ 
	such that $|S| < c(G-S)-1$. 
	Let $|S| = s$, $c(G-S) = c$, which implies $s \leq c-2$. 
	For positive integers $n_1 \geq n_2 \geq \cdots \geq n_{s+2} \geq 1$ 
	with $\sum_{i=1}^{s+2} n_i = n-s$, 
	$G$ is a spanning subgraph of 
	$\hat{G} =K_s \vee (K_{n_1} \cup K_{n_2} \cup \ldots \cup K_{n_{s+2}})$. 
	By Lemma \ref{le1}, we have
	\begin{equation}\label{eq1}
		\lambda_1(D(\hat{G})) \leq \lambda_1(D(G))
	\end{equation}
	with equality if and only if $G \cong \hat{G}$. 
	Let $\tilde{G} = K_s \vee (K_{n-2s-1} \cup (s+1)K_1)$. 
	By Lemma \ref{le3}, we have
	\begin{equation}\label{eq2}
	\lambda_1(D(\tilde{G})) \leq \lambda_1(D(\hat{G}))
	\end{equation}
	with equality if and only if $\hat{G} \cong \tilde{G}$.
	
	we divide the proof into the following three cases.
	
	\vspace{1.8mm}
	\noindent\textbf{Case 1.} $s = \delta$
	\vspace{1mm}
	
 In this case, $\tilde{G} = K_\delta \vee (K_{n-2\delta-1} \cup (\delta+1)K_1)=G^*$. 
	From (\ref{eq1}) and (\ref{eq2}), 
	$\lambda_1(D(G^*)) \leq \lambda_1(D(G))$ with equality if and only if $G \cong G^*$. 
 Since $\lambda_1(D(G)) \leq \lambda_1(D(G^*))$, 
		we have $\lambda_1(D(G)) = \lambda_1(D(G^*))$. 
	Therefore, $G \cong G^*=K_\delta \vee (K_{n-2\delta-1} \cup (\delta+1)K_1)$. 
	Let $S = V(K_\delta)$, we can calculate
	
	\[
	\frac{|S|}{c(G-S)-1} = \frac{\delta}{\delta+1} < 1.
	\]
	Thus $\tau(G) < 1$, and $G$ is not a $1$-tough graph.
	
	\vspace{1.8mm}
	\noindent\textbf{Case 2.} $s \geq \delta+1$
	\vspace{1mm}
	
	The distance matrix \(D(G^*)\) of $G^*$ is
	\[
	\bordermatrix{
		& \delta+1    & n-2\delta-1  & \delta \cr
		\hfill \delta+1     & 2(J-I)    &  2J        & J      \cr
		\hfill n-2\delta-1  & 2J        & J-I        & J      \cr
		\hfill \delta     & J         & J          & J-I    \cr
	}.
	\]
	
	We can partition the vertex set of $G^*$ as $V(G^*)=V((\delta+1)K_1) \cup V(K_{n-2\delta-1}) \cup V(K_\delta)$. The quotient matrix is
	\[
	R_\delta =
	\begin{pmatrix}
		2\delta & 2(n-2\delta-1) & \delta \\
		2(\delta+1) & n-2\delta-2 & \delta \\
		\delta+1 & n-2\delta-1 & \delta-1
	\end{pmatrix}.
	\]
 By a direct calculation, the characteristic polynomial of $R_\delta$ is
	\begin{equation}
		\begin{aligned}\label{eq3}
			P(R_{\delta},x)=&x^{3}-(n+\delta-3)x^{2}-((2\delta+5)n-5\delta^{2}-6(\delta+1))x\\
			&+(\delta^{2}-\delta-4)n-2\delta^{2}(\delta-1)+6\delta+4.
		\end{aligned}
	\end{equation}
Clearly, the vertex partition of $G^*$ is equitable. 
	By Lemma \ref{le2}, $\lambda_1(D(G^*)) = \lambda_1(R_\delta)$ is the largest root of $P(R_\delta,x) = 0$. 
 Also, note that $D(\tilde{G})$ has the equitable quotient matrix $R_s$, 
		which is obtained by replacing $\delta$ with $s$ in $R_\delta$. 
		Thus, by (\ref{eq3}), we have 
		\begin{align*}
			P(R_{s},x)=&x^{3}-(n+s-3)x^{2}-((2s+5)n-5s^{2}-6(s+1))x\\
			&+(s^{2}-s-4)n-2s^{2}(s-1)+6s+4.
		\end{align*}
		According to Lemma~\ref{le2}, the largest root of 
		$P(R_{s},x)=0$ equals $\lambda_1(D(\tilde{G}))$. 
		We are to verify $P(R_{\delta},x)-P(R_{s},x)>0$. By a simple computation, we obtain
		\begin{align*}
P(R_{\delta},x)-P(R_{s},x)=(s-\delta)f(x),
		\end{align*}
		where $f(x)=x^{2}+(2n-5\delta-5s-6)x
		-(s+\delta-1)n+2(s^{2}+s\delta+\delta^{2}-s-\delta)-6$. 
		Then the symmetry axis of 
		$f(x)$ is $x=-n+\frac{5}{2}s+\frac{5}{2}\delta + 3$, which implies that 
		$f(x)$ is increasing in the interval $[-n+\frac{5}{2}s+\frac{5}{2}\delta + 3,+\infty)$. 
		
		By Lemma~\ref{le3} and \(n\geq9\delta\), we note that 
		\begin{align*}
			\lambda_1(D(G^{*}))&\geq\frac{2W(G^{*})}{n}=\frac{(\delta + 1)(2n - \delta - 2)+(n - 2\delta - 1)(n + \delta)+\delta(n - 1)}{n}\\
			&=(n+\delta + 1)+\delta-\frac{3\delta^{2}+5\delta+ 2}{n}\geq n+\delta + 1
		\end{align*}
		Since $\delta+1\leq s$ and $n\geq s+c\geq 2s+2$, we have $\delta+1\leq s\leq\frac{n-2}{2}$ and 
		$$
		-n+\frac{5}{2}s+\frac{5}{2}\delta + 3\leq(n+\delta + 1)-\frac{3}{2}s+\frac{3}{2}\delta - 2<n+\delta + 1,
		$$
		which implies that 
		\begin{align*}
			f(x)&\geq f(n+\delta + 1)\\
			&=2s^{2}-(6n + 3\delta + 7)s + 3n^{2}-2\delta n - n - 2\delta^{2}-11\delta - 11\\
			&\geq 2\left(\frac{n - 2}{2}\right)^{2}-(6n + 3\delta + 7)\left(\frac{n - 2}{2}\right)+3n^{2}-2\delta n - n - 2\delta^{2}-11\delta - 11\\
			&=\frac{1}{2}n^{2}-(\frac{7}{2}\delta+\frac{1}{2})n - 2\delta^{2}-8\delta - 2\\
			&\geq2\delta^{2}-3\delta - 1\\
			&>0.
		\end{align*}
		It follows that $P(R_{\delta},x)>P(R_{s},x)$ for $x\geq n+\delta + 1$. 
		Recall that $\lambda_1(D(G^*))> n+\delta+1$ and 
		$\lambda_1(D(G^*))$ is the largest root of $P(R_\delta,x) = 0$. 
		Thus, $\lambda_1(D(G^*)) < \lambda_1(D(\tilde{G}))$. By (\ref{eq1}) and (\ref{eq2}), we have
	\[
	\lambda_1(D(G^*)) < \lambda_1(D(\tilde{G})) \leq \lambda_1(D(\hat{G})) \leq \lambda_1(D(G)),
	\]
	a contraction.
	
	\vspace{1.8mm}
	\noindent\textbf{Case 3.} $1 \leq s \leq \delta-1$
	\vspace{1mm}
	
	Recall that for positive integers $n_1 \geq n_2 \geq \cdots \geq n_{s+2} \geq 1$ with $\sum_{i=1}^{s+2} n_i = n-s$, $G$ is a spanning subgraph of $\hat{G} = K_s \vee (K_{n_1} \cup K_{n_2} \cup \ldots \cup K_{n_{s+2}})$. Note that $\delta(\hat{G}) \geq \delta(G) = \delta$, so $n_{s+2}-1+s \geq \delta$, implying $n_1 \geq n_2 \geq \cdots \geq n_{s+2} \geq \delta-s+1$. 
	We claim that $n_1 \geq 2(\delta-s+1)$. 
	Suppose to the contrary that $n_1 \leq 2\delta-2s+1$. 
	Notice that $n_1 \geq n_2 \geq \cdots \geq n_{s+2} \geq 1$ 
	and $1 \leq s \leq \delta-1$. Then we have
	\begin{align*}
		n&=s + n_{1}+n_{2}+\cdots+n_{s + 1}+n_{s + 2}\\
		&\leq s+(s + 2)(2\delta - 2s + 1)\\
		&=-2s^{2}+(2\delta - 2)s + 4\delta + 2:=n(s).
	\end{align*}
	The symmetry axis of $n(s)$ is $s_{\text{axis}} = \frac{\delta-1}{2}$. 
	When $\delta \geq 3$, we have $1 \leq s_{\text{axis}} \leq \delta-1$ and 
	$$
	n(s)\leq - 2\left(\frac{\delta - 1}{2}\right)^2
	+(2\delta - 2)\left(\frac{\delta - 1}{2}\right)+4\delta + 2
	=\frac{1}{2}\delta^2+3\delta+\frac{5}{2}.
	$$
	When $2 \leq \delta \leq 3$, we have $s_{\text{axis}} \leq 1$ and
	$$
	n(s)\leq - 2 + 2\delta-2 + 4\delta + 2 <6\delta.
	$$
	Both cases contradict $n \geq \max\{9\delta, \frac{\delta^2}{2} + 3\delta + 3\}$. Therefore, the assumption is false, and we conclude that $n_1 \geq 2(\delta-s+1)$.
	
	Let $G' = K_s \vee (K_{n-s-(\delta-s+1)(s+1)} \cup (s+1)K_{\delta-s+1})$. By Lemma \ref{le4}, we have
	\begin{equation}\label{eq4}
		\lambda_1(D(G')) \leq \lambda_1(D(\hat{G}))
	\end{equation}
	with equality if and only if $G' \cong \hat{G}$.
	
	\vspace{1.8mm}
	\noindent\textbf{Subcase 1.} $s=1.$
	\vspace{1mm}
	
	In this case, $G' = K_1 \vee (K_{n-2\delta-1} \cup 2K_\delta)$. Its distance matrix \(D(G')\) is
	\[
	\bordermatrix{
		& \delta    & \delta  & n-2\delta-1 &1 \cr
		\hfill \delta 	&J-I & 2J & 2J & J \cr
		\hfill \delta &2J & J-I & 2J & J \cr
		\hfill n-2\delta-1 &2J & 2J & J-I & J \cr
		\hfill 1 &J & J & J & O
	}.
	\]
	
	Recall that $G^* = K_\delta \vee (K_{n-2\delta-1} \cup (\delta+1)K_1)$. Let $\mathbf{x}$ be the Perron vector of $D(G^*)$. By symmetry, $\mathbf{x}$ takes the same value within each of the vertex sets $V((\delta+1)K_1)$, $V(K_{n-2\delta-1})$, and $V(K_\delta)$. Denote the components of $\mathbf{x}$ corresponding to these vertex sets as $x_1$, $x_2$, and $x_3$, respectively. From $D(G^*)\mathbf{x} = \lambda_1(D(G^*))\mathbf{x}$, we have
	\begin{equation*}
		\begin{aligned}
			\lambda_1(D(G^*))x_1 &= 2\delta x_1 + 2(n-2\delta-1)x_2 + \delta x_3;\\
			\lambda_1(D(G^*))x_2 &= 2(\delta+1)x_1 + (n-2\delta-2)x_2 + \delta x_3;\\ 
			\lambda_1(D(G^*))x_3 &= (\delta+1)x_1 + (n-2\delta-1)x_2 + (\delta-1)x_3.
		\end{aligned}
	\end{equation*}
	Thus,
	\begin{equation*}
		\begin{aligned}
		 \lambda_1(D(G^*))(2x_3 - x_1) &= 2x_1 + (\delta-2)x_3;\\
		 \lambda_1(D(G^*))(x_2 - x_3) &= (\delta+1)x_1 - (x_2 - x_3).
		\end{aligned}
	\end{equation*}
	Since $\delta \geq 2$, we have $2x_1 + (\delta-2)x_3 \geq 2x_1 > 0$ and $(\lambda_1(D(G^*)) + 1)(x_2 - x_3) = (\delta+1)x_1 > 0$. Therefore, $2x_3 \geq x_1 > 0$ and $x_2 > x_3 > 0$.
	
	By simple calculation, we obtain that $D(G') - D(G^*)$ is
	\[ 
	\bordermatrix{
		& \delta  &1  & \delta-1  & n-2\delta-1 &1 \cr
		\hfill \delta	&-(J-I) & O & J & O & O \cr
		\hfill 1	&O & O & O & O & O \cr
		\hfill \delta-1	&J & O & O & J & O \cr
		\hfill n-2\delta-1	&O & O & J & O & O \cr
		\hfill 1	&O & O & O & O & O
	}. \]
	From this, we obtain that 
	\begin{align*}
		&\lambda_1(D(G')) - \lambda_1(D(G^{*})) \geq \mathbf{x}^{T}(D(G') - D(G^{*}))\mathbf{x}\\
		=&-\delta(\delta - 1)x_1^{2}+2\delta(\delta - 1)x_1x_3+(\delta - 1)(n - 2\delta - 1)x_2x_3+(\delta - 1)(n - 2\delta - 1)x_2^{2}\\
		=&\delta(\delta - 1)(2x_3 - x_1)x_1+(\delta - 1)(n - 2\delta - 1)x_2x_3+(\delta - 1)(n - 2\delta - 1)x_2^{2}\\
		>& 0.
	\end{align*}
	Therefore, $\lambda_1(D(G')) > \lambda_1(D(G^*))$.
	
	\vspace{1.8mm}
	\noindent{\bf Subcase 2.} $2 \leq s \leq \delta-1$
	\vspace{1mm}
	
	In this case, $G' = K_s \vee (K_{n-s-(\delta-s+1)(s+1)} \cup (s+1)K_{\delta-s+1})$. 
	We can partition the vertex set of $G'$ as $V(G') = V((s+1)K_{\delta-s+1}) \cup V(K_{n-s-(\delta-s+1)(s+1)}) \cup V(K_s)$. The quotient matrix is
	\[ R_{s,\delta} = 
	\begin{pmatrix}
		(\delta-s) + 2s(\delta-s+1) & 2[n-s-(\delta-s+1)(s+1)] & s \\
		2(s+1)(\delta-s+1) & n-s-(\delta-s+1)(s+1)-1 & s \\
		(s+1)(\delta-s+1) & n-s-(\delta-s+1)(s+1) & s-1
	\end{pmatrix}. \]
	After calculation, the characteristic polynomial of $R_{s,\delta}$ is
	\begin{equation}
		\begin{aligned}\label{eq5}
			P(R_{s,\delta},x)=&x^{3}+[s^{2}-(\delta + 1)s - n + 3]x^{2}+[2s^{4}-(4\delta + 2)s^{3}\\&+(2\delta^{2}-3\delta + 2n - 3)s^{2}
			+(5\delta^{2}-2n\delta + n + 5\delta)s + 3\delta^{2}-3n\delta + 6\delta - 5n + 6]x\\ &- s^{5}+(2\delta + 3)s^{4}
			-(\delta^{2}+3\delta + n)s^{3}+(n\delta - 6\delta + 2n - 5)s^{2}\\
			&+(4\delta^{2}-n\delta + 4\delta + 2n)s+ 3\delta^{2}-3n\delta + 6\delta - 4n.
		\end{aligned}
	\end{equation}
	Clearly, the vertex partition of $V(G')$ is equitable. By Lemma \ref{le2}, $\lambda_1(D(G')) = \lambda_1(R_{s,\delta})$ is the largest root of $P(R_{s,\delta},x) = 0$. 
	Note that $G^*$ contains the proper subgraph $K_{n-\delta-1}$. 
	Then $\lambda_1(D(G^*)) > \lambda_1(D(K_{n-\delta-1})) = n - \delta - 2$. Combining (\ref{eq3}) and (\ref{eq5}), we have
	\begin{align*}
		&P(R_{\delta},n - \delta - 2)-P(R_{s,\delta},n - \delta - 2)\\
		=&(\delta - s)[3sn^{2}+(4\delta + 2s^{3}-(2\delta + 3)s^{2}-(9\delta + 9)s + 4)n - s^{4}-(\delta + 1)s^{3}\\
		&+(2\delta^{2}+6\delta + 4)s^{2}+(6\delta^{2}+11\delta + 5)s - 5\delta^{2}-9\delta - 4]\\
		:=&(\delta - s)g(n).
	\end{align*}
	Since $2 \leq s \leq \delta-1$ and $\delta \geq s + 1 \geq 3$, the symmetry axis of $g(n)$ is
	$$
	\frac{-2s^{3}+(2\delta + 3)s^{2}+(9\delta + 9)s - 4\delta - 4}{6s}
	\leq-\frac{1}{3}s^{2}+\frac{2\delta + 3}{6}s+\frac{3\delta + 3}{2}
	:=\varphi(s).
	$$
	The symmetry axis of $\varphi(s)$ is  $s_{\text{axis}} = \frac{2\delta + 3}{4}$. 
Therefore, 
		$$
		\varphi(s)\leq \varphi\left(\frac{2\delta + 3}{4}\right)
		=\frac{1}{12}\delta^2 + \frac{7}{4}\delta + \frac{27}{16}
		< \frac{1}{2}\delta^2 + 2\delta + 2.
		$$
	This implies that $g(n)$ is monotonically increasing for 
	$n \in [\frac{1}{2}\delta^2 + 2\delta + 2, +\infty)$. Since $2 \leq s$ and $\delta \geq s + 1 \geq 3$, we have
	\begin{align*}
		g(n) \geq& g\left(\frac{1}{2}\delta^2 + 2\delta + 2\right) \\
		=& \frac{\delta}{4}[\delta(3s\delta^2 - (4s^2 - 6s - 8)\delta + 4s^3 - 14s^2 + 6s + 20) \\
		& + 12s^3 - 16s^2 - 4s + 28] - (s-1)(s^3 - 2s^2 + 1) + 3 \\
		\geq& \frac{\delta}{4}[\delta(3s^3 - 6s^2 + 23s + 28) + 12s^3 - 16s^2 - 4s + 28] - (s-1)(s^3 - 2s^2 + 1) + 3 \\
		\geq& \frac{\delta}{4}(3s^4 + 6s^3 + 7s^2 + 24s + 28) - (s-1)(s^3 - 2s^2 + 1) + 3 \\
		\geq& \frac{1}{4}(3s^5 + 2s^4 + 19s^3 + 16s^2 + 24s + 16) \\
		>& 0.
	\end{align*}
	It follows that 
	\begin{equation}\label{eq6}
		P(R_{\delta},n - \delta - 2) > P(R_{s,\delta},n - \delta - 2).
	\end{equation}
	
	On the other hand, for $x \in [n - \delta - 2, +\infty)$ and $s \geq 2$, we have
	\begin{align*}
		P'(R_{\delta},x) - P'(R_{s,\delta},x) &= (\delta - s)[(2s - 2)x + 2s^3 - (2\delta + 2)s^2 \\
		&\quad - (5\delta - 2n + 3)s + 2\delta + n] \\
		&\geq (\delta - s)[(2s - 2)(n - \delta - 2) + 2s^3 - (2\delta + 2)s^2 \\
		&\quad - (5\delta - 2n + 3)s + 2\delta + n] \\
		&= (\delta - s)[2s^3 - (2\delta + 2)s^2 + (4n - 7\delta - 7)s + 4\delta - n + 4] \\
		&:= (\delta - s)h(s).
	\end{align*}
	Then we prove that $h(s) > 0$ for $2 \leq s \leq \delta - 1$. Note that 
	$h'(s) = 6s^2 - 4(\delta + 1)s + 4n - 7\delta - 7$. 
	The symmetry axis of $h'(s)$ is $s_{\text{axis}} = \frac{\delta + 1}{3}$. 
	When $3 \leq \delta \leq 5$, we have $s_{\text{axis}} \leq 2$ and  
	\begin{align*}
		h'(s) &\geq h'(2) = 4n - 15\delta + 9 \geq 2\delta^2 - 3\delta + 3 > 0.
	\end{align*}
	When $\delta > 5$, we have $2 < s_{\text{axis}} \leq \delta - 1$ and 
	\begin{align*}
		h'(s) &\geq h'\left(\frac{\delta + 1}{3}\right) = 4n - \frac{2}{3}\delta^2 - \frac{25}{3}\delta - \frac{23}{3} \geq \frac{4}{3}\delta^2 + 3\delta + \frac{7}{3} > 0.
	\end{align*}
	This implies that $h(s)$ is monotonically increasing for $2 \leq s \leq \delta - 1$. Combining with $n \geq 9\delta$ and $\delta \geq 3$,
	we have
	\begin{align*}
		h(s) &\geq h(2) = 7n - 18\delta - 2 \geq 45\delta - 2 > 0.
	\end{align*}
	It follows that 
	\begin{equation}\label{eq7}
		P'(R_{\delta},x) > P'(R_{s,\delta},x).
	\end{equation}
	
	Next, we consider
	\(
	P'(R_{\delta},x) = 3x^2 - 2(\delta + n - 3)x + 5\delta^2 + (6 - 2n)\delta - 5n + 6.
	\)
	Note that $n \geq 9\delta$ and $\delta \geq 3$. 
	The symmetry axis for $P'(R_{\delta},x)$ is
$$
		\frac{\delta + n - 3}{3}\leq (n - \delta - 2) - \frac{14}{3}\delta + 1 
		< n - \delta - 2.
		$$
	Thus,
	\begin{equation}\label{eq8}
		\begin{aligned}
			P'(R_{\delta},x) &\geq P'(R_{\delta},n - \delta - 2) = n^2 - (8\delta + 7)n + 10\delta^2 + 16\delta + 6 \\
			&\geq 19\delta^2 - 47\delta + 6> 0.
		\end{aligned}
	\end{equation}
	Therefore, $P(R_{\delta},x)$ is monotonically increasing for $x \in [n - \delta - 2, +\infty)$. 
	
	Combining (\ref{eq6}), (\ref{eq7}), and (\ref{eq8}), we conclude that $\lambda_1(D(G')) > \lambda_1(D(G^*))$.
	Combining (\ref{eq1}), (\ref{eq4}), Subcase 1, and Subcase 2, we have
	\begin{align*}
		\lambda_1(D(G^*)) < \lambda_1(D(G')) \leq \lambda_1(D(\hat{G})) \leq \lambda_1(D(G)),
	\end{align*}
	{a contraction. This completes the proof of Theorem~\ref{main1}.}
	\hspace*{\fill}$\Box$

	\section{Proof of Theorem \ref{main2}} 
	In this section, we give the proof of Theorem \ref{main2}.
	\begin{proof}
	Let $G^*=K_{\tau-1} \vee (K_{n-\tau} \cup K_1)$.
		Suppose to the contrary that $G$ is not a $\tau$-tough graph. 
		Then there exists a non-empty vertex subset $S \subset V(G)$ 
		such that $|S| < \tau(c(G-S)-1)$. Let $|S| = s$ and $c(G-S) = c$. 
		Since $\tau \geq 2$, we have $s \leq \tau(c-1)-1$. 
		For positive integers $n_1 \geq n_2 \geq \cdots \geq n_c \geq 1$ 
		with $\sum_{i=1}^c n_i = n-\tau(c-1)+1$, 
		$G$ is a spanning subgraph of 
		$\hat{G} = K_{\tau(c-1)-1} \vee (K_{n_1} \cup K_{n_2} 
		\cup \ldots \cup K_{n_c})$. By Lemma \ref{le1}, we have
		\begin{equation}\label{eq9}
			\lambda_1(D(\hat{G})) \leq \lambda_1(D(G))
		\end{equation}
		with equality if and only if $G \cong \hat{G}$.
		
		Let $\tilde{G} = K_{\tau(c-1)-1} \vee (K_{n-(\tau+1)(c-1)+1} \cup (c-1)K_1)$. By Lemma \ref{le3}, we have
		\begin{equation}\label{eq10}
			\lambda_1(D(\tilde{G})) \leq \lambda_1(D(\hat{G}))
		\end{equation}
		with equality if and only if $\hat{G} \cong \tilde{G}$. 
		
		We consider two cases based on the range of $c$.

		\vspace{1.8mm}
		\noindent\textbf{Case 1.} $c = 2$
		\vspace{1mm}

		In this case, $\tilde{G}=K_{\tau-1}\vee (K_{n-\tau}\cup K_1)=G^*$. 
		From (\ref{eq9}) and (\ref{eq10}), we have 
		$\lambda_1(D(G^*)) \leq \lambda_1(D(G))$ 
		with equality if and only if $G \cong G^*$. 
	Since  $\lambda_1(D(G)) \leq \lambda_1(D(G^*))$, we have
		$\lambda_1(D(G^*)) = \lambda_1(D(G))$. 
	This implies that
		$G \cong G^*=K_{\tau-1} \vee (K_{n-\tau} \cup K_1)$. 
		Let $S = V(K_{\tau-1})$, then we can calculate
		\[
		\frac{|S|}{c(G-S)-1} = \tau - 1 < \tau.
		\]
		Thus $\tau(G) < \tau$, and $G \cong K_{\tau-1} \vee (K_{n-\tau} \cup K_1)$ is not a $\tau$-tough graph.
		
		\vspace{1.8mm}
		\noindent\textbf{Case 2.} $c \geq 3$
		\vspace{1mm}
		
		Consider $\tilde{G} = K_{\tau(c-1)-1} \vee (K_{n-(\tau+1)(c-1)+1} \cup (c-1)K_1)$. 
	By Lemma \ref{le5}, we have
			$$
			\lambda_1(D(\tilde{G}))\geq\frac{2W(\tilde{G})}{n}
			=\frac{-(2\tau + 1)(c-1)^2 + (2n+1)(c-1) + n^2 - n}{n}
			:=\frac{\varPhi(c)}{n}.
			$$
			From $\sum_{i=1}^c n_i = n-\tau(c-1)+1$, we have 
			$n-\tau(c-1)+1 \geq c$. Then $2 \leq c-1 \leq \frac{n}{\tau+1}$. 
			Let $x = c-1$. Thus $2 \leq x \leq \frac{n}{\tau+1}$ and 
			\[
			\Phi(x) = -(2\tau + 1)x^2 + (2n+1)x + n^2 - n.
			\]
			Note that $n \geq 4\tau^2 + 5\tau + 1$. Therefore, 
			$$
			\varPhi\left(\frac{n}{\tau + 1}\right)-\varPhi(2)
			=\frac{(n-2\tau-2)(n-4\tau^{2}-5\tau-1)}{(\tau + 1)^2}
			\geq 0.
			$$
		Thus, $\min \Phi(x) = \Phi(2)$. We can get
		\begin{equation}\label{e-1}
			\lambda_1(D(\tilde{G})) \geq \frac{\Phi(2)}{n}
			= \frac{n^2+3n-8\tau-2}{n}
			\geq n+3.
		\end{equation}
		
		Now consider $G^* = K_{\tau-1} \vee (K_{n-\tau} \cup K_1)$. Its distance matrix \(D(G^*)\) is
		\[
		\bordermatrix{
			& \tau-1    & n-\tau& 1\cr
			\hfill \tau-1 		&J-I & J & J \cr
			\hfill n-\tau 		&J & J-I & 2J \cr
			\hfill 1 		&J & 2J & O\cr
		}.
		\]
		By direct calculation, we have
		\begin{align*}
			W^{(2)}(G^{*})&=\sum_{1\leq i<j\leq n}d_{ij}^2(G^{*})\\
			&=\frac{(\tau - 1)(n - 1)+(n - \tau)(n+2)+(\tau - 1)+4(n - \tau)}{2}\\
			&=\frac{1}{2}n^2+\frac{5}{2}n - 3\tau
		\end{align*}
		By Lemma \ref{le6} and $\tau \geq 2$, we have
		\begin{align*}
			\lambda_1(D(G^{*}))&\leq\sqrt{\frac{2(n - 1)W^{(2)}(G^{*})}{n}}=\sqrt{\frac{(n - 1)(n^2 + 5n - 6\tau)}{n}}\\
			&<\sqrt{\frac{n^3 + 4n^2 - 5n}{n}}=\sqrt{(n + 2)^2 - 9}\\
			&\leq n + 2
		\end{align*}
		Therefore, by (\ref{e-1}), $\lambda_1(D(G^*)) \leq n+2 < \lambda_1(D(\tilde{G}))$. Combining (\ref{eq9}) and (\ref{eq10}), we have
		\[
		\lambda_1(D(G^*)) < \lambda_1(D(\tilde{G})) \leq \lambda_1(D(\hat{G})) \leq \lambda_1(D(G)),
		\]
		a contraction. This completes the proof of Theorem~\ref{main2}
	\end{proof}
	
	\section{Proof of Theorem \ref{main3}}
	In this section, we give the proof of Theorem \ref{main3}
	\begin{proof}
		Let $G^*=K_1 \vee (K_{n-\frac{1}{\tau}-2} \cup (\frac{1}{\tau}+1)K_1)$. 
		Suppose to the contrary that $G$ is not a $\tau$-tough graph. 
		Then there exists a non-empty vertex subset $S \subset V(G)$ such that $|S| < \tau(c(G-S)-1)$. Let $|S| = s$ and $c(G-S) = c$. It follows that $c \geq \frac{s}{\tau} + 2$. 
		For positive integers $n_1 \geq n_2 \geq \cdots \geq n_{\frac{s}{\tau}+2} \geq 1$ with $n_1 + n_2 + \cdots + n_{\frac{s}{\tau}+2} = n - s$, 
		$G$ is a spanning subgraph of 
		$\hat{G} = K_s \vee (K_{n_1} \cup K_{n_2} \cup \cdots \cup K_{n_{\frac{s}{\tau}+2}})$. By Lemma \ref{le1}, we have
		\begin{equation}\label{eq11}
			\lambda_1(D(\hat{G})) \leq \lambda_1(D(G)).
		\end{equation}
		with equality if and only if $G \cong \hat{G}$.
		
		Let $\tilde{G} = K_s \vee (K_{n-s-(\frac{s}{\tau}+1)} \cup (\frac{s}{\tau}+1)K_1)$. By Lemma \ref{le3}, we have
		\begin{equation}\label{eq12}
			\lambda_1(D(\tilde{G})) \leq \lambda_1(D(\hat{G})).
		\end{equation}
		with equality if and only if $\hat{G} \cong \tilde{G}$.
		
		Since $S \neq \varnothing$, we have $s \geq 1$. We consider two cases based on the range of $s$.
		
		\vspace{1.8mm}
		\noindent\textbf{Case 1.} $s = 1.$
		\vspace{1mm}
		
		In this case, $\tilde{G}=K_1 \vee (K_{n-1-(\frac{1}{\tau}+1)} \cup (\frac{1}{\tau}+1)K_1)=G^*$. 
		From (\ref{eq11}) and (\ref{eq12}), we have 
		$\lambda_1(D(G^*)) \leq \lambda_1(D(G))$ 
		with equality if and only if \(G \cong G^*\). 
		Since $\lambda_1(D(G)) \leq \lambda_1(D(G^*))$,  
		we obtain $\lambda_1(D(G^*)) = \lambda_1(D(G))$. 
		This implies that 
		$G \cong G^*=K_1 \vee (K_{n-1-(\frac{1}{\tau}+1)} \cup (\frac{1}{\tau}+1)K_1)$. 
		Let $S = V(K_1)$. Then we can calculate
		\[
		\frac{|S|}{c(G-S)-1} = \frac{1}{1 + \frac{1}{\tau} + 1 - 1} < \tau.
		\]
		Thus $\tau(G) < \tau$, 
		and $G \cong K_1 \vee (K_{n-1-(\frac{1}{\tau}+1)} \cup (\frac{1}{\tau}+1)K_1)$ is not a $\tau$-tough graph.
		
		\vspace{1.8mm}
		\noindent\textbf{Case 2.} $s \geq 2.$
		\vspace{1mm}
		
		Let $\tilde{G} = K_s \vee (K_{n-s-(\frac{s}{\tau}+1)} \cup (\frac{s}{\tau}+1)K_1)$. 
		We can partition the vertex set of $\tilde{G}$ as $V(\tilde{G}) = V((\frac{s}{\tau}+1)K_1) \cup V(K_{n-s-(\frac{s}{\tau}+1)}) \cup V(K_s)$. The quotient matrix is
		\[
		R_{\tau,s} = 
		\begin{pmatrix}
			2\frac{s}{\tau} & 2(n-s-(\frac{s}{\tau}+1)) & s \\
			2(\frac{s}{\tau}+1) & n-s-(\frac{s}{\tau}+1)-1 & s \\
			\frac{s}{\tau}+1 & n-s-(\frac{s}{\tau}+1) & s-1
		\end{pmatrix}.
		\]
		The characteristic polynomial of $R_{\tau,s}$ is
		\begin{align*}
			P(R_{\tau,s},x) &= x^3 + \frac{-n\tau^2 - s\tau + 3\tau^2}{\tau^2}x^2+ \frac{-2ns\tau - 5n\tau^2 + 3s^2\tau + 2s^2 + 3s\tau + 3s\tau^2 + 6\tau^2}{\tau^2}x \\
			&\quad + \frac{ns^2\tau + ns\tau^2 - 2ns\tau - 4n\tau^2 - s^3\tau - s^3 - \tau^2s^2 + s^2\tau + 2s^2 + 2s\tau^2 + 4s\tau + 4\tau^2}{\tau^2}.
		\end{align*}
		Clearly, the vertex partition of $\tilde{G}$ is equitable. 
		By Lemma \ref{le2}, $\lambda_1(D(\tilde{G})) = \lambda_1(R_{\tau,s})$ is the largest root of $P(R_{\tau,s},x) = 0$. 
		Also, note that $D(G^*)$ has the equitable quotient matrix $R_\tau=R_{\tau,1}$, 
		which is obtained by replacing $s$ with $1$ in $R_{\tau,s}$. 
		According to Lemma~\ref{le2}, the largest root of 
		$P(R_\tau,x)=0$ equals $\lambda_1(D(G^*))$. 
		We are to verify $P(R_{\tau},x)-P(R_{\tau,s},x)>0$. By a simple computation, we obtain
		\begin{align*}
			P(R_{\tau},x)-P(R_{\tau,s},x)=\frac{s - 1}{\tau^{2}}\alpha(x).
		\end{align*}
		where $\alpha(x)=(\tau x^{2}+(2n\tau - 3s\tau - 2s - 3\tau^{2}-6\tau - 2)x- n(s\tau + \tau^{2}-\tau)+s^{2}\tau + s^{2}+s\tau^{2}-s - \tau^{2}-4\tau - 1)$. 
		Then the symmetry axis of 
		$\alpha(x)$ is $x=-n+\frac{3}{2}s+\frac{s}{\tau}+\frac{3\tau}{2}+3+\frac{1}{\tau}$, which implies that 
		$\alpha(x)$ is increasing in the interval $[-n+\frac{3}{2}s+\frac{s}{\tau}+\frac{3\tau}{2}+3+\frac{1}{\tau},+\infty)$. Consider $G^* = K_1 \vee (K_{n-\frac{1}{\tau}-2} \cup (\frac{1}{\tau}+1)K_1)$. 
		By Lemma \ref{le5} and $n \geq 4\tau + \frac{1}{\tau} + 5$, we have
		\[
		\begin{aligned}
			\lambda_{1}(D(G^{*}))&\geq\frac{2W(G^{*})}{n}=\frac{(1+\frac{1}{\tau})(2n-3)+(n-\frac{1}{\tau}-2)(n+\frac{1}{\tau})+n - 1}{n}\\
			&=n+\frac{1}{\tau}+1+\frac{1}{\tau}-\frac{4}{n}-\frac{5}{n\tau}-\frac{1}{n\tau^{2}}\\
			&\geq n+\frac{1}{\tau}+1
		\end{aligned}
		\]

			Since  $\sum_{i=1}^{\frac{s}{\tau}+2} n_i = n - s$, we have $n \geq s+\frac{s}{\tau}+2$. Then $2 \leq s \leq \frac{n-2}{1+\frac{1}{\tau}}$ and 
			\[
			-n+\frac{3}{2}s+\frac{s}{\tau}+\frac{3\tau}{2}+3+\frac{1}{\tau}
			\leq(n+\frac{1}{\tau}+1)-2(s+\frac{s}{\tau}+2)+\frac{3}{2}s+\frac{s}{\tau}+\frac{3}{2}\tau+2
			\leq n+\frac{1}{\tau}+1
			\]
			which implies that 
			\begin{align*}
				\alpha(x) \geq \alpha\left(n + \frac{1}{\tau} + 1\right) 
				&= (\tau + 1)x^2 + \left(-4n\tau - 2n + \tau^2 - 3\tau - 6 - \frac{2}{\tau}\right)x \\
				&\quad - 4\tau^2 - 12\tau - 7 - \frac{1}{\tau} + 3n^2\tau + n(-4\tau^2 - \tau + 2) \\
				&\geq (\tau + 1)\left(\frac{n-2}{1+\frac{1}{\tau}}\right)^2 + \frac{n-2}{1+\frac{1}{\tau}}\left(-4n\tau - 2n + \tau^2 - 3\tau - 6 - \frac{2}{\tau}\right) \\
				&\quad -4\tau^2 - 12\tau - 7 - \frac{1}{\tau} + 3n^2\tau + n(-4\tau^2 - \tau + 2) \\
				&\geq \frac{1}{\tau^2 + \tau}\left(n^2\tau^2 - (3\tau^4 + 4\tau^3 + \tau^2)n - 6\tau^4 - 6\tau^3 - 7\tau^2 - 4\tau - 1\right) \\
				&\geq \frac{1}{(\tau^2 + \tau)(\tau + 1)}\left(-12\tau^6 - 33\tau^5 - 14\tau^4 + 24\tau^3 + 22\tau^2 + 5\tau\right)\\
				&>0,
			\end{align*}
		It follows that $P(R_{\tau},x)>P(R_{\tau,s},x)$ for $x\geq n+\frac{1}{\tau}+1$. 
	Observed that $\lambda_1(D(G^*))> n+\frac{1}{\tau}+1$ and 
	$\lambda_1(D(G^*))$ is the largest root of $P(R_{\tau},x) = 0$. 
	Thus, $\lambda_1(D(G^*)) < \lambda_1(D(\tilde{G}))$. Combining (\ref{eq11}) and (\ref{eq12}), we have
		\[
		\lambda_1(D(G^*)) < \lambda_1(D(\tilde{G})) \leq \lambda_1(D(\hat{G})) \leq \lambda_1(D(G))
		\]
		a contraction. This completes the proof of Theorem~\ref{main3}.
	\end{proof}

	\vspace{5mm}
	\noindent
	\vspace{3mm}
	

	

\end{document}